\documentclass{amsart}

\usepackage[dvips]{graphicx}
\usepackage{amsthm}
\usepackage{amsmath}
\usepackage{amsfonts}
\usepackage{amssymb}
\usepackage{epic}
\usepackage{graphicx}
\usepackage{epstopdf}
\usepackage{tabularx}
\usepackage{array}
\usepackage{multirow}
\usepackage{longtable}

\newcolumntype{L}[1]{>{\raggedright\let\newline\\\arraybackslash\hspace{0pt}}m{#1}}
\newcolumntype{C}[1]{>{\centering\let\newline\\\arraybackslash\hspace{0pt}}m{#1}}
\newcolumntype{R}[1]{>{\raggedleft\let\newline\\\arraybackslash\hspace{0pt}}m{#1}}

\newtheorem{definition}{Definition}

\DeclareMathOperator{\Av}{Av}
\DeclareMathOperator{\st}{st}
\DeclareMathOperator{\maj}{maj}
\DeclareMathOperator{\inv}{inv}

\newcommand{\ignore}[1]{}

\begin{document} 

\title{Inversion Polynomials for Permutations Avoiding Consecutive Patterns} 
\date{}


\author{Naiomi T. Cameron \& Kendra Killpatrick}
\address{Lewis \& Clark College}
\email{ncameron@lclark.edu}

\begin{abstract}
	
In 2012, Sagan and Savage introduced the notion of $st$-Wilf equivalence for a statistic $st$ and for sets of permutations that avoid particular permutation patterns which can be extended to generalized permutation patterns.  In this paper we consider $\inv$-Wilf equivalence on sets of two or more consecutive permutation patterns.  We say that two sets of generalized permutation patterns $\Pi$ and $\Pi'$ are $\inv$-Wilf equivalent if the generating function for the inversion statistic on the permutations that simultaneously avoid all elements of $\Pi$ is equal to the generating function for the inversion statistic on the permutations that simultaneously avoid all elements of $\Pi'$.  

In 2013, Cameron and Killpatrick gave the inversion generating function for Fibonacci tableaux which are in one-to-one correspondence with the set of permutations that simultaneously avoid the consecutive patterns $321$ and $312.$  In this paper, we use the language of Fibonacci tableaux to study the inversion generating functions for permutations that avoid $\Pi$ where $\Pi$ is a set of five or fewer consecutive permutation patterns.  In addition, we introduce the more general notion of a strip tableaux which are a useful combinatorial object for studying consecutive pattern avoidance.  We go on to give the inversion generating functions for all but one of the cases where $\Pi$ is a subset of three consecutive permutation patterns and we give several results for $\Pi$ a subset of two consecutive permutation patterns.  

\end{abstract}

\maketitle

\section{Introduction}

For any $\pi \in S_n$ and $\sigma \in S_k$, we say that $\pi$ {\it {contains a copy of $\sigma$}} if $\pi$ has a subsequence that is order isomorphic to $\sigma$.  If $\pi$ contains no subsequence order isomorphic to $\sigma$ then we say that $\pi$ {\it {avoids}} $\sigma$.  This is the notion of classical pattern avoidance and we write $\sigma=\sigma_1\sigma_2\cdots\sigma_k$ with dashes between each of the elements $\sigma_1$, $\sigma_2$, $\dots$, $\sigma_k$ to indicate that these elements need not be adjacent in $\pi$.  For the sake of clarity in this paper, we will refer to a pattern of the sort $\sigma_1 - \sigma_2 - \sigma_3 - \cdots - \sigma_k$ as a {\em standard pattern} of length $k$.  

Generalized permutation patterns were first introduced by Babson and Steingrimmson in 2000 \cite{Bst}.  In a {\em generalized permutation pattern} $\sigma$, we may require some of the letters in $\sigma$ to be adjacent, which will be indicated by the lack of a dash between these letters in $\sigma$.  For example, a $3-1-2$ pattern in $\pi$ would be any subsequence $\pi_i$, $\pi_j$, $\pi_k$ for $i<j<k$ with $\pi_j < \pi_k < \pi_i$, while a $31-2$ pattern in $\pi$ would be any subsequence $\pi_i$, $\pi_{i+1}$, $\pi_j$ for $i+1 < j$ with $\pi_{i+1} < \pi_j < \pi_i$.  A pattern with no internal dashes will be referred to as a {\em consecutive pattern}.  A number of interesting results on generalized permutation patterns were obtained by Claesson \cite{Cla}, including relations to several well studied combinatorial structures such as set partitions, Dyck paths, Motzkin paths and involutions.  

The study of permutations that avoid a given pattern has been a topic of much interest in the past decade.  Early work on generalized permutation patterns was done by Elizalde and Noy \cite{ElN} and Kitaev \cite{Kit} in 2003.  Elizalde and Noy gave the first systematic study of permutations that avoid generalized permutation patterns of length 3, primarily focusing on consecutive patterns of length 3.  Kitaev considered simultaneous avoidance of two or more consecutive patterns of length 3.  He gave enumerative results for the number of permutations that avoid any subset of four such consecutive patterns and for all but one subset of three such consecutive patterns.  He also gave several results for subsets of size two.  Kitaev and Mansour \cite{KiM} continued this work, enumerating the permutations that avoid the remaining size three subset, and Aldred, Atkinson and McCaughan \cite{AAM} completed the enumeration of the permutations that avoid the subsets of size two.  In 2012, Elizalde and Noy \cite{ElN2} studied pattern avoidance for subsets of consecutive patterns of length greater than three.

In 2011, Dokos et al.\ \cite{DDJ} began to examine the generating functions for permutation statistics on permutations that avoid standard patterns of length three.  Their work focuses on two well-known permutation statistics, the inversion statistic and the major index, giving generating functions and recursions for these statistics on certain sets of pattern-avoiding permutations and leaving open the recursions for other sets.  In 2012, Cheng et al. \cite{CEK} continued this work by giving a recursion for the inversion polynomial on permutations avoiding the standard pattern $3-2-1.$  

In 2010, Cameron and Killpatrick \cite{CaK} gave a recursion for the inversion polynomial on Fibonacci tableaux, which are in one-to-one correspondence with the set of permutations that avoid both the consecutive patterns $321$ and $312.$  In this paper, we use Fibonacci tableaux and a more general notion of strip tableaux to determine the inversion polynomials on permutations that avoid multiple consecutive patterns of length three.  In Section \ref{background}, we give all the necessary background and definitions needed to study this type of multi-avoidance.  In particular, we define a strip tableau which is related to Fibonacci tableaux and turns out to be a useful combinatorial object in the study of multi-avoidance of generalized patterns.  Also in Section \ref{background} we make a connection to the $q$-Catalan polynomials.  A subset of Fibonacci tableaux with all columns of height two is counted by the Catalan numbers and there is a nice relation between the inversion polynomial on this subset and the well known $q$-Catalan polynomials.  In Section \ref{invpoly}, we provide combinatorial proofs for the inversion polynomials of all sets of permutations avoiding three or more consecutive patterns.  We also provide some results for the inversion polynomial on permutations avoiding two or fewer consecutive patterns.

\section{Background}\label{background}

\vspace{.1in}

\subsection{Wilf Equivalence}

Let $\Pi$ be a set of standard permutation patterns in $S_k$ where $k\leq n$ and define $\Av_n(\Pi)$ as the set of permutations in $S_n$ which avoid every pattern in $\Pi$.  Two sets of permutation patterns $\Pi$ and $\Pi'$ are said to be {\it {Wilf equivalent}} if $\vert \Av_n(\Pi)\vert = \vert \Av_n(\Pi')\vert$.  If $\Pi$ and $\Pi'$ are Wilf equivalent, we write $\Pi \equiv \Pi'$.  It is well-known that for any two standard permutation patterns $\tau, \sigma$ of length three, $\vert \Av_n(\tau)\vert = \vert \Av_n(\sigma) \vert$.  That is, there is just one Wilf equivalence class on the set of standard permutation patterns of length 3.

Now we extend this notion to the case where $\Pi$ and $\Pi'$ are sets of generalized permutation patterns.  As above, we will say $\Pi$ and $\Pi'$ are {\it {Wilf equivalent}} if $\vert \Av_n(\Pi) \vert = \vert \Av_n(\Pi') \vert$. 

Throughout this paper we will utilize some basic operations on permutations, namely the {\it {inverse}}, the {\it {reverse}} and the {\it {complement}}.  For a permutation $\pi = \pi_1 \pi_2 \cdots \pi_n$, the inverse is the standard group-theoretic inverse operation on permutations, the reverse is 
\[
\begin{array}{cccccc}
\pi^r& =& \pi_n& \cdots& \pi_2 &\pi_1
\end{array}
\]
and the complement is 
\[
\begin{array}{cccccc}
\pi^c &= &(n+1)-a_1 & (n+1)-a_2 & \cdots & (n+1)-a_n.
\end{array}\]
 
Since each of these operations gives a bijection from $S_n$ to itself, it is easy to see that for any standard permutation pattern $\sigma$, $\vert \Av_n(\sigma) \vert = \vert \Av_n(\sigma^{-1}) \vert = \vert \Av_n(\sigma^r) \vert = \vert \Av_n(\sigma^c) \vert$.  For consecutive patterns $\sigma$ we have $\vert \Av_n(\sigma) \vert = \vert \Av_n(\sigma^r) \vert = \vert \Av_n(\sigma^c) \vert$.  
   
\subsection{st-Wilf Equivalence}
Sagan and Savage \cite{SSa} recently defined a $q$-analogue of Wilf equivalence by considering any permutation statistic $\st$ from $\uplus_{n \geq 0} S_n \rightarrow \mathbb{N}$, where $\mathbb{N}$ is the set of nonnegative integers, and letting 
\[
F_n^{\st}(\Pi;q) = \sum_{\sigma \in \Av_n(\Pi)} q^{\st(\sigma)}.
\]
For $\Pi$ and $\Pi'$ subsets of permutations, they defined $\Pi$ and $\Pi'$ to be {\it {$\st$-Wilf equivalent}} if $F_n^{\st}(\Pi;q) = F_n^{\st}(\Pi';q)$ for all $n \geq 0$.  In this case, we write $\Pi \stackrel{\st}{\equiv} \Pi'$. We will use $[\Pi]_{\st}$ to denote the $\st$-Wilf equivalence class of $\Pi$.  If we set $q=1$ in the generating function above we have $F_n^{\st}(\Pi;1) = \vert \Av_n(\Pi) \vert$, thus $\st$-Wilf equivalence implies Wilf equivalence.  While Sagan and Savage defined this notion for $\Pi$ and $\Pi'$ subsets of standard permutation patterns, we can use the same definition for $\Pi$ and $\Pi'$ sets of generalized permutation patterns.  

In \cite{DDJ}, Dokos et al.\ give a thorough investigation of $\st$-Wilf equivalence for standard permutation patterns of length 3 for both the major index, {\it {$\maj$}}, and the inversion statistic, {\it {$\inv$}}.  Through a relatively straightforward map on permutations that takes the major index to another well known Mahonian statistic, the charge statistic, one can give a similarly thorough investigation for the charge statistic and Killpatrick makes this explicit in \cite{Kil}.  The main focus of this paper will be the inversion statistic on permutations avoiding sets of consecutive patterns of length three.

Given a permutation $\sigma  = \sigma_1 \sigma_2 \cdots \sigma_n \in S_n$ where $\sigma_i=\sigma(i)$, we define an {\bf inversion} to be a pair $(i, j)$ such that $i<j$ and $\sigma_i > \sigma_j$.  Then the {\bf inversion statistic}, $\inv(\sigma)$, is the total number of inversions in $\sigma$, i.e.
\[
\inv(\sigma) = \sum_{\substack{i<j \\ \sigma_i>\sigma_j}} 1
\]
  For example, for $\sigma = \begin{array}{ccccccccc} 3&2&8&5&7&4&6&1&9 \end{array}$, $\inv(\sigma) = 15$. 

It is relatively straightforward to prove that for $\sigma \in S_n$, 
\begin{equation}
\inv(\sigma) = \inv(\sigma^{-1}) = \binom{n}{2}-\inv(\sigma^r) = \binom{n}{2}-\inv(\sigma^c).
\label{invoperations}
\end{equation}

 We use $I_n(\Pi;q)$ to denote the inversion polynomial on the set of permutations that avoid the set $\Pi$ of generalized permutation patterns.  That is,
\[
I_n(\Pi;q)=F_n^{\inv}(\Pi;q) = \sum_{\sigma \in \Av_n(\Pi)} q^{\inv(\sigma)}
\]

\subsection{Fibonacci tableaux}\label{Fibtableaux}

Fibonacci tableaux arise from the Fibonacci lattice first defined by Stanley \cite{St1} in 1975.  We define a Fibonacci shape to be a set of consecutive columns of heights one or two.  The Ferrers diagram for a Fibonacci shape is formed by replacing a $1$ with a single dot and a $2$ with two dots.  The size of the Fibonacci shape is the sum of the $1$'s and $2$'s that make up the shape.  For example, the Ferrers diagram for the Fibonacci shape $\mu = 122121$ looks like
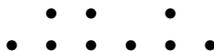
\[
\begin{array}{cccccc}
 &\bullet&\bullet& &\bullet& \\
\bullet&\bullet&\bullet&\bullet&\bullet&\bullet
\end{array}
\]
and has size 9.

\begin{definition}
Given a Fibonacci shape $\mu$ of size $n$ and its corresponding Ferrers diagram, a {\bf {standard Fibonacci tableau}} is a filling of the Ferrers diagram with the integers $1$ through $n$ such that the bottom row decreases from left to right and every column decreases from bottom to top. 
\end{definition}

An example of a standard Fibonacci tableau of shape $\mu = 122121$ is
\[
\begin{array}{cccccc}
 &3&4& &2& \\
9&8&7&6&5&1
\end{array}
\] 

\begin{definition}
We define the {\bf {column-reading word}}, $w_c(T)$, for a standard Fibonacci tableau $T$ by reading the columns from bottom to top, right to left.  
\end{definition}
For the tableau above, $w_c(T) = 152674839$.   
We define $\inv(T)$ for $T$ a Fibonacci tableau to be $\inv(w_c(T))$.

Due to the restrictions on the standard Fibonacci tableaux (i.e. columns must decrease from bottom to top and the bottom row must decrease from left to right), we obtain an immediate bijection between the set of column reading words for standard Fibonacci tableaux and the set of permutations which avoid both the generalized patterns $312$ and $321$.  In \cite{CaK}, the authors used standard Fibonacci tableaux to prove the following recurrence for the inversion polynomial on permutations that avoid $\Pi = \{321, 312\}$:
\begin{equation}
I_n(\Pi;q) = I_{n-1}(\Pi;q) + (q + q^2 + \cdots + q^{n-1})I_{n-2}(\Pi;q).
\label{class17}
\end{equation}

Now we extend the notion of a Fibonacci tableau to a more general strip tableau which will allow us to model more sets of generalized pattern avoiding permutations.  To begin, we must define a general strip shape.

\begin{definition}  A {\bf {strip shape of size $n$}} is a set of n contiguous cells such that if each cell were numbered 1 through $n$ consecutively beginning with the leftmost and uppermost cell, then cell $i$ must be below or to the right of cell $i-1$ for all $i$.  
\end{definition}

For example,
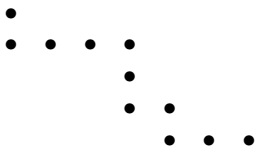
\[
\begin{array}{ccccccc}
\bullet& & & & & & \\
\bullet&\bullet&\bullet&\bullet& & & \\
 & & &\bullet& & & \\
 & & &\bullet&\bullet& & \\
 & & & &\bullet&\bullet&\bullet
\end{array}\]

is a strip shape of size 11.

\begin{definition}  A {\bf {standard strip tableau of size $n$}} is a filling of a strip shape of size $n$ with the numbers 1 through $n$ such that each column decreases from bottom to top and each row decreases from left to right.
\end{definition}

For example,
\[
T=\begin{array}{ccccccc}
4& & & & & & \\
8&6&5&1& & & \\
 & & &9& & & \\
 & & &10&7& & \\
 & & & &11&3&2
\end{array}\]

is a standard strip tableau of size 11.

\begin{definition}
The {\bf {column-reading word}} $w_c(T)$ of a standard strip tableau $T$ is obtained by reading the columns of $T$ from bottom to top, right to left.  
\end{definition}
For the tableau $T$ above, $w_c(T) = \begin{array}{ccccccccccc} 2&3&11&7&10&9&1&5&6&8&4\end{array}$.   

Note that the map $w_c$ which takes a standard strip tableau $T$ of size $n$ to $w_c(T)$ is injective and therefore reversible.  The inverse map $w_c^{-1}$ takes a permutation $\pi \in S_n$ to its corresponding strip tableau $T$ in the following way.  Reading $\pi=\pi_1\pi_2\cdots\pi_n$ from left to right, place $\pi_1$ in a cell and for $i=2,\dots,n$, if $\pi_i<\pi_{i-1}$, place $\pi_i$ on top of $\pi_{i-1}$; otherwise place $\pi_i$ to the left of $\pi_{i-1}$.  Observe that if $i$ is a descent of $\pi$, then $\pi_i$ would appear in a column of height greater than one in the corresponding strip tableau $T_\pi$.

\begin{figure}
\begin{center}
\begin{tabularx}{\textwidth}{|p{1cm}|p{1cm}|X|}
\hline\hline
Class & $\Pi$ & Conditions on the set of strip tableaux $T$ corresponding to $Av_n(\Pi)$ \\\hline\hline
I &321 & every column of $T$ has at most 2 entries\\\hline
II &312 &  for every $i\in T$ which appears at the rightmost end of a row, the entry below $i$ is less than the entry to the left of $i$ (alternatively, the two entries adjacent to $i$ decrease left to right OR ``right corners decrease") \\\hline
III &231 &  for every $i\in T$ which appears at the leftmost end of a row, the two entries adjacent to $i$ decrease left to right (``left corners decrease")\\\hline
IV &213 & for every $i\in T$ which appears at the rightmost end of a row, the two entries adjacent to $i$ increase left to right (``right corners increase")\\\hline
V &132 & for every $i\in T$ which appears at the leftmost end of a row, the two entries adjacent to $i$ increase left to right (``left corners increase")\\\hline
VI &123 &  every row of $T$ has at most 2 entries\\\hline\hline
\end{tabularx}
\end{center}
\caption{The set of strip tableaux of size $n$ corresponding to the generalized pattern avoiding permutations in $Av_n(\Pi)$.}
\label{restrictedstriptableaux}
\end{figure}

If $\Pi$ is a nonempty set of generalized patterns, then the set $\Av_n(\Pi)$ is in one-to-one correspondence with a proper subset of the set of strip tableaux of size $n.$  The table in Figure \ref{restrictedstriptableaux} describes the restrictions on strip tableaux that correspond to $\Av_n(\Pi)$ when $\Pi$ consists of a single consecutive pattern of length three.  As noted in section \ref{Fibtableaux}, there is a bijection between $\Av_n(321,312)$ and the set of standard Fibonacci tableaux of size $n$.  Hence, we have a bijection between the set of standard Fibonacci tableaux and the set of strip tableaux in both Class I and Class II of Figure \ref{restrictedstriptableaux}.  This bijection can be seen immediately by example.  If a strip tableaux $T$ meets the criteria of both Class I and Class II, then to obtain the corresponding Fibonacci tableaux, simply drop each column of $T$ to the bottom level while maintaining its row position.  The result will be a tableau in which each column has height at most 2, all columns decrease from bottom to top and the bottom row decreases from left to right.

In our work with the inversion polynomial for certain subsets of strip tableaux, we will need to consider various cases where the strip tableau begins with a (leftmost) column of height $k$ for some $k$ or ends with a (rightmost) column of height $k$, leading us to the following definition.

\begin{definition}
Let $I_n^k(\Pi;q)$ denote the inversion polynomial for the subset of strip tableaux that avoid $\Pi$ and that begin with a column of height $k$.  Similarly, let ${\widetilde{I_n^{k}}}(\Pi;q)$ denote the inversion polynomial for the subset of strip tableaux that avoid $\Pi$ and end with a column of height $k$. 
\label{defOandTw}
\end{definition}
\subsection{$q$-Catalan polynomials}

The well known Catalan numbers are defined as 
\[
C_n = \frac{1}{n+1} \binom{2n}{n}.
\]

In addition to having an explicit formula, the Catalan numbers are
known to satisfy the recurrence \[ C_n=\sum_{i=1}^{n} C_{i-1} C_{n-i}.
\]

A {\em Dyck path\/} is a lattice path in $\mathbb{Z}^2$ from $(0,0)$ to
$(n,n)$ consisting of only steps in the positive $x$ direction (EAST
steps) and steps in the positive $y$ direction (NORTH steps) such that
there are no points $(x, y)$ on the path for which $x > y$.  Let $D_n$
denote the set of Dyck paths from $(0,0)$ to $(n,n)$.  The Catalan number $C_n$ is known to count the number of Dyck paths
from $(0,0)$ to $(n,n)$, thus $C_3 = 5$.  The {\em length\/} of a Dyck
path is the number of NORTH steps in the path, thus a Dyck path $\pi
\in D_n$ has length $n$.

Given a Dyck path $\pi \in D_n$ the
{\em area\/} statistic, $a(\pi)$, is the number of squares that lie
below the path and completely above the diagonal.  The generating function for the area statistic on Dyck paths $\pi \in
D_n$, \[ \sum_{\pi \in D_n} q^{a(\pi)} = C_n(q) \] is called the
$q$-Catalan polynomial and was first defined by Carlitz and Riordan
\cite{CR}.  Specializing $q=1$ in the $q$-Catalan polynomial gives the
usual Catalan number $C_n$.  It is known that \[ C_n(q) =
\sum_{i=1}^{n} q^{i-1}C_{i-1}(q) C_{n-i}(q).  \]

There is an explicit bijection between Dyck paths and standard two row tableaux of size $2n$.  If one labels the steps of the Dyck path from $1$ to $2n$ beginning with the lower left corner and ending at the upper right corner, then one can form a standard two row tableaux of size $2n$ by putting the labels of the NORTH steps in the top row and the EAST steps in the bottom row (both increasing from left to right).

Such tableaux are in one-to-one correspondence with the number of standard Fibonacci tableaux that have all columns of height two and whose elements in the top row of each column also decrease from left to right.  The inversion polynomial for such tableaux is exactly $q^k$ times the $q$-Catalan polynomial $C_k(q)$.

To prove this, note that the standard Fibonacci tableau of size $2k$ with all columns of height two having the smallest inversion statistic is the tableau 
\[
T = 
\begin{array}{cccc} 
2k-1 & 2k-3 & \cdots& 1 \\
2k & 2k-2 & \cdots & 2
\end{array}\]
which corresponds to the Dyck path that alternates NORTH and EAST steps.  This tableau has an inversion statistic of $k$ while the Dyck path has an area statistic of 0.  To form a new Dyck path, change any adjacent EAST-NORTH pair of steps to a NORTH-EAST pair of steps.  This increases the area statistic on the Dyck path by 1.  In the corresponding Fibonacci tableau, this operation is equivalent to changing the position of two adjacent numbers $i$ and $i+1$, which maintains the fact that the tableau is a standard Fibonacci tableau with elements in the top row decreasing from left to right.  Swapping adjacent numbers creates one new inversion in the column reading word of the tableau, thus the inversion statistic increases by 1 as well.



\section{Inversion Polynomials for Generalized Pattern Avoiding Permutations}\label{invpoly}

In this section, we will give formulas for $I_n(\Pi;q)$ where $\Pi$ contains two or more consecutive patterns of length three.  

First, we consider permutations that avoid all but one consecutive pattern of length 3.  For $n=3$ there is exactly one such permutation and the inversion statistic can be calculated directly, so we assume $n\geq 4$. There is exactly one such permutation, the all increasing permutation, that avoids all consecutive patterns of length three except 123.  This permutation has an inversion statistic of 0.  There is also only one such permutation, the all decreasing permutation, that avoids all consecutive patterns of length three except 321.  This permutation has an inversion statistic of $q^{\binom{n}{2}}$.  There are no other permutations that avoid five consecutive patterns of length three.

We now consider $\Pi$ where $\Pi$ contains four or fewer consecutive patterns of length three.  For each Wilf equivalence class below, we provide an argument for the inversion polynomial of the indicated representative $\Pi$.  As in \cite{DDJ}, to obtain inversion polynomials for the remaining members of the Wilf class we rely on equation (\ref{invoperations}).  The inversion polynomials $I_n(\Pi;q)$ are recorded in Table \ref{bigtable}.

\subsection*{Class 1 - $\Pi=\{321,312, 132, 123\}$}\label{Class1}

There are only two permutations in $\Av_n(\Pi)$.  If $n$ is odd they are
\[
\begin{array}{cccccccc}
\frac{n+1}{2}&\frac{n+3}{2} &\frac{n-1}{2} & \frac{n+5}{2}  & \frac{n-3}{2}  & \cdots& n & 1 
\end{array}
\]
and 
\[
\begin{array}{cccccccc}
\frac{n+1}{2}&\frac{n-1}{2} &\frac{n+3}{2} & \frac{n-3}{2}  & \frac{n+5}{2}  & \cdots & 1 & n
\end{array}
\]
If $n$ is even they are
\[
\begin{array}{ccccccc}
\frac{n}{2}&\frac{n+2}{2} &\frac{n-2}{2} & \frac{n+4}{2}  & \frac{n-4}{2}  & \cdots & n 
\end{array}
\]
and
\[
\begin{array}{ccccccc}
\frac{n+2}{2}&\frac{n}{2}&\frac{n+4}{2} & \frac{n-2}{2}  & \frac{n+6}{2} & \cdots & 1.
\end{array}
\]

One can compute the inversion polynomials directly for these permutations, thus 
\begin{equation}
I_n(\Pi;q)=\begin{cases}q^{\left(\frac{n+1}{2}\right)\left(\frac{n-1}{2}\right)} + q^{\left(\frac{n-1}{2}\right)\left(\frac{n-1}{2}\right)}&\text{if $n$ is odd}\cr
q^{\left(\frac{n-2}{2}\right)\left(\frac{n}{2}\right)} + q^{\left(\frac{n}{2}\right)\left(\frac{n}{2}\right)}&\text{if $n$ is even.}
\end{cases}
\label{class1}
\end{equation}

Using the language of Fibonacci tableaux, since the permutations that avoid this set of patterns avoid both $321$ and $312$ they correspond to a subset of Fibonacci tableaux.  Since permutations in this set also avoid the $123$ pattern, a Fibonacci tableau in this subset must also have at most two columns of height one which occur (if they do occur) as the first and/or last column in the tableau.  Permutations that avoid this subset also avoid the $132$ pattern which gives the condition that the elements in the top of each column of the Fibonacci tableau must increase from left to right.  

Thus if $n$ is odd, the two Fibonacci tableaux are
\[
\begin{array}{ccccc}
1&2&\dots&\frac{n-1}{2}& \\
n&n-1& &\frac{n+1}{2}+1&\frac{n+1}{2}
\end{array}\]
and
\[
\begin{array}{ccccc}
 &1&2&\dots&\frac{n-1}{2}\\
n&n-1&n-2& &\frac{n+1}{2}
\end{array}.\]

If $n$ is even, the two Fibonacci tableaux are
\[
\begin{array}{cccc}
1&2&\dots&\frac{n}{2}\\
n&n-1& &\frac{n}{2}+1
\end{array}
\]
and
\[
\begin{array}{cccccc}
 &1&2&\dots&\frac{n}{2}-1& \\
n&n-1&n-2& &\frac{n}{2}+1&\frac{n}{2}
\end{array}\]
One can clearly see that the reading words of these tableaux give the permutations above.

\subsection*{Class 2 - $\Pi=\{321,312, 231,132\}$}\label{Class2}

Again for this class, there are only two permutations that avoid all four of these patterns and they are
\[
\begin{array}{ccccc}
1&2&\cdots&n-1&n
\end{array}\] and
\[
\begin{array}{ccccccc}
2&1&3&4&\cdots&n-1&n
\end{array}.
\]
Hence,
\begin{equation}
I_n(\Pi;q)=1+q.
\label{class2}
\end{equation}

In terms of Fibonacci tableaux, permutations that avoid $\Pi$ are in correspondence with Fibonacci tableaux that have at most one column of height two which must appear in the rightmost position of the tableau.  This corresponds to the following Fibonacci tableaux:

\[
\begin{array}{ccccc}
n&n-1&\cdots&2&1
\end{array}\]
\[
\begin{array}{cccccc}
 & & & & &1\\
n&n-1&n-2&\cdots&3&2
\end{array}
\]

\subsection*{Class 3 -$\Pi=\{312, 213, 231, 132\}$}\label{Class3}

The only two permutations to avoid all four permutations in this subset are the all increasing or all decreasing permutation thus we have
\begin{equation}
I_n(\Pi;q)=1+q^{\binom{n}{2}}.
\label{class3}
\end{equation}

Since permutations that avoid this set of patterns do not necessarily avoid the 321 pattern they do not correspond with a subset of Fibonacci tableaux.  Instead we must use the more general language of strip tableaux to find the appropriate combinatorial interpretation.  Since these permutations either strictly increase or strictly decrease, they correspond to the two strip tableaux that consist of either a single row or a single column, respectively.

\subsection*{Class 4 - $\Pi=\{321,312, 213,123\}$}
If $n=3$ there are only two permutations in $\Av_3(\Pi)$ and they are 132 and 231, thus 
\begin{equation}
I_3(\Pi;q) = q + q^2.
\label{class4}
\end{equation}
The corresponding Fibonacci tableaux are $\begin{array}{cc} 2& \\3&1 \end{array}$ and $\begin{array}{cc} 1& \\3&2 \end{array}$.  If $n > 3$, there are no permutations in $\Av_n(\Pi)$.

\subsection*{Class 5 - $\Pi=\{321,312, 213, 132\}$}

Since the permutations in $\Av_n(\Pi)$ must avoid $321, 312$ and $213$, this set of permutations is in correspondence with standard Fibonacci tableaux that have at most one column of height two which, if it exists, must be the first column in the Fibonacci tableau.

If the Fibonacci tableau has all columns of height one then it corresponds with the all increasing permutation and has an inversion statistic of 0.  If the Fibonacci tableau begins with a column of height two then any of the numbers $1$ through $n-2$ can occur as the topmost element in this column (not $n-1$ since the corresponding permutation must also avoid 132).  If $k$ is the number in the top of the first column of height two then the corresponding permutation has an inversion statistic of $n-k$ thus
\begin{equation}
I_n(\Pi;q)=1+q^2+\cdots+q^{n-1}=\frac{1-q+q^2-q^n}{1-q}.
\label{class5}
\end{equation}

\subsection*{Class 6 - $\Pi=\{321,312, 123,231\}$}\label{Class6}

Since the permutations in $\Av_n(\Pi)$ must avoid $321, 312$ and $123$, this set of permutations is in one to one correspondence with standard Fibonacci tableaux that have either no columns of height one, exactly one column of height one that occurs in the rightmost or leftmost position, or exactly two columns of height one that occur in the rightmost and leftmost positions.  Thus, if $n$ is even, we have standard Fibonacci tableaux with all columns of height two or with two columns of height one in the rightmost and the leftmost column positions.  If $n$ is odd, we have standard Fibonacci tableaux with a single column of height one in either the rightmost position or the leftmost position.  In addition, the condition that the permutations in $\Av_n(\Pi)$ must avoid $231$ implies that the elements in the top of each column of the tableaux must decrease from left to right.  

Given that \begin{equation} I_n(\Pi;q) = I_n^{1}(\Pi;q) + I_n^{2}(\Pi;q)\label{OplusTw}\end{equation} we consider the case of $n$ odd and $n$ even separately.

\begin{enumerate}
\item[Case 1.]  Suppose $n=2k+1$ for some $k$.  If a permutation in $\Av_n(\Pi)$ corresponds to a Fibonacci tableau with a column of height one in the rightmost position, then by the above restrictions on the standard Fibonacci tableaux, the $1$ in the tableau must occur in this column.  The set of standard Fibonacci tableaux obtained by removing this rightmost column of height one and relabeling the elements in the resulting tableau with the numbers $1$ through $n-1$ (in the same relative order as they appeared in the original tableau) are in bijection with the set of standard Fibonacci tableaux of size $2k$ with the restriction that all columns have height two and the elements in the top of each column decrease from left to right.  (Note:  these tableaux are counted by the Catalan number $C_k$.).  For the original tableau, the 1 is in the first position of the column reading word so this 1 does not contribute to the inversion statistic and removing it and relabeling preserves the inversion statistic.  Thus $I_{2k+1}^{2}(\Pi;q) = I_{2k}^{2}(\Pi;q)=q^k C_k(q)$.

If the Fibonacci tableau has a column of height one in the leftmost position, then by the above restrictions on the standard Fibonacci tableaux, the $n$ in the tableau must occur in this column.  The set of standard Fibonacci tableaux obtained by removing this leftmost column of height one are in bijection with the set of standard Fibonacci tableaux of size $2k$ with the restriction that all columns have height 2 and the elements in the top of each column decrease from left to right.  (Note:  these tableaux are again counted by the Catalan number $C_k$.).  In the original tableau, the $n$ is in the last position when reading the word of the standard tableau, thus this $n$ does not contribute to the inversion statistic and removing this $n$ preserves the inversion statistic.  Thus $I_{2k+1}^1(\Pi;q) = I_{2k}^{2}(\Pi;q)=q^k C_k(q)$ and by (\ref{OplusTw}), we have
\begin{equation}
I_{2k+1}(\Pi;q)= 2 I_{2k}^{2}(\Pi;q)= 2q^kC_k(q).
\label{class6odd}
\end{equation}

\item[Case 2.] 

Suppose $n=2k$ for some $k$.  If the Fibonacci tableau has columns of height one in both the leftmost and the rightmost columns, then these columns must contain $n$ and $1$ respectively.  The set of standard Fibonacci tableaux obtained by removing these two columns of height one and relabeling the elements in the resulting tableau with the numbers $1$ through $n-2$ (in the same relative order as they appeared in the original tableau) are in bijection with the set of standard Fibonacci tableaux of size $2k-2$ with the restriction that all columns have height two and the topmost row decreases from left to right.  In the original tableau, the $n$ is in the last position and the $1$ is in the first position when reading the word of the standard tableau, thus the $1$ and the $n$ do not contribute to the inversion statistic and removing these two elements and relabeling preserves the inversion statistic.  Thus 
\[
I_{2k}^1(\Pi;q) = I_{2k-2}^{2}(\Pi;q)=q^{k-1}C_{k-1}(q).
\]

If the Fibonacci tableau has only columns of height two, then let $m$ be the smallest integer such that the first $m$ columns of the tableau contain the elements $n$, $n-1$, $\cdots$, $n-2m+1$.  In these first $m$ columns, the $n$ is in the bottom cell of the leftmost column and the $n-2m+1$ is in the upper cell of the rightmost column.  Remove these two cells and slide the remaining cells in the second row one column to the right to obtain a new tableau with $k-1$ columns containing the numbers $n-1$, $n-2$, $\cdots$, $n-2m+2$ such that the top row decreases from left to right.  With proper relabeling, one can see that these tableaux are in one-to-one correspondence with the 2-row tableaux counted by the Catalan number $C_{m-1}$.  The change in the inversion statistic is $m$.  Thus 

\begin{align*}
I_{2k}^{2}(\Pi;q) &= \sum_{m=1}^k q^m I_{2m-2}^{2}(\Pi;q) I_{2k-2m}^{2}(\Pi;q)\\
&=\sum_{m=1}^k q^m q^{m-1}C_{m-1}(\Pi;q) q^{k-m}C_{k-m}(\Pi;q)
\end{align*}

 and by (\ref{OplusTw}), we have
\begin{align}
I_{2k}(\Pi;q) &= I_{2k-2}^{2}(\Pi;q) + \sum_{m=1}^k q^m I_{2m-2}^{2}(\Pi;q) I_{2k-2m}^{2}(\Pi;q)\\
& = q^{k-1}C_{k-1}(q)  + \sum_{m=1}^k q^m q^{m-1}C_{m-1}(q) q^{k-m}C_{k-m}(q)\\
&= q^{k-1}C_{k-1}(q)+q^kC_k(q)
\label{class6even}
\end{align}
\end{enumerate}

\subsection*{Class 7 - $\Pi= \{321,312, 231\}$}\label{Class7}

Since the permutations in $\Av_n(\Pi)$ must avoid $321$, $312,$ and $231,$ this set of permutations is in one-to-one correspondence with standard Fibonacci tableaux for which the elements in the top of each column must decrease from left to right.   (Note: the element in a column of height one is both in the top row of its column and in the bottom row of its column.) 

If the Fibonacci tableau begins with a column of height one, then there must be an $n$ in this column and $n$ appears as the last element in the word of the tableau.  Removing this $n$ will not change the inversion statistic for this permutation and will give a tableau of size $n-1$ with the given restrictions.  Therefore, $$I_n^1(\Pi;q) = I_{n-1}(\Pi;q).$$

If the Fibonacci tableau begins with a column of height two, let $k$ be the smallest integer such that the first $k-1$ columns have height two and the $k$th column has height one.  Since the elements in the top row of each column and the bottom row of each column must decrease from left to right, the numbers $n$ through $n-2(k-1)+1$ must be in the first $k-1$ columns and $n-2(k-1)$ must be in column $k$ of height 1.  Thus the tableaux formed from the remaining columns to the right of column $k$ corresponds to a permutation of size $n-2(k-1)-1 = n-2k + 1$.  Since $n-2k+2$ is in column $k$, removal of this element does not change the inversion statistic and also gives a tableau in the first $k-1$ columns that corresponds (with relabeling) to a tableau of size $2(k-1)$ with all columns of height two whose elements in the top row of each column decrease from left to right.  These tableaux are counted by the Catalan number $C_{k-1}$.  This gives us the inversion polynomial $$q^{k-1}C_{k-1}(\Pi;q) I_{n-(2k-1)}(\Pi;q),$$  and we have 
$$I_{n}^{2}(\Pi;q)=\sum_{k=1}^{\lfloor n/2 \rfloor}{q^{k-1} C_{k-1}(q)\cdot I_{n-2k+1}(\Pi;q)}.$$ Thus, with initial conditions $I_2^1(\Pi;q)=1$ and $I_2^{2}(\Pi;q)=I_3^{2}(\Pi;q)=q$, we have

\begin{eqnarray}
I_n(\Pi;q)&=&I_{n-1}(\Pi;q) + \sum_{k=1}^{\lfloor n/2 \rfloor}{q^{k-1} C_{k-1}(q)\cdot I_{n-2k+1}(\Pi;q)}
\label{class7}
\end{eqnarray}

\subsection*{Class 8 - $\Pi=\{321,312, 213\}$}\label{Class8}

Since the permutations in $\Av_n(\Pi)$ can have at most one descent in the last position, the corresponding Fibonacci tableaux contain either $n$ columns of height one or an initial column of height two followed by $n-2$ columns of height one.  The permutation corresponding to the tableau with $n$ columns of height one is $\pi=123\cdots n$ which has an inversion statistic of 0.  There are $n-1$ Fibonacci tableaux of shape $211\cdots1$ and these can be formed by placing $n$ in the bottom cell of the column of height two, then choosing a $k$, $1 \leq k \leq n-1$, to be the element in the top cell of the column of height two, then placing the remaining elements in the tableau in decreasing order from left to right.  If $k$ is the element in the top cell of the column of height two, then the corresponding permutation has $n-k$ inversions.  Thus 
\begin{equation}
I_n(\Pi;q)=1+q+\cdots+q^{n-1}=\frac{1-q^n}{1-q}.
\label{class8}
\end{equation}

\subsection*{Class 9 - $\Pi=\{312, 231, 132\}$}\label{Class9}

Since permutations in $\Av_n(\Pi)$ avoid both the $231$ and $132$ patterns, once there is an ascent in the permutation then there can be no further descents.  This corresponds to a strip tableau with shape a single row followed by a single column.  Since the permutations also avoid $312$, the single right corner (if it exists in the shape) must decrease.  

Since ${\widetilde{I_n^k}}(\Pi; q)$ is the inversion generating function for all permutations of length $n$ that avoid $\Pi$ and have a last column of height $k$, we have $I_n(\Pi;q) = \sum_{k=1}^n {\widetilde{I_n^k}}(\Pi;q)$.  

The permutation corresponding to $k=1$ is $\begin{array}{cccc} 1&2&\cdots&n \end{array}$ and has an inversion statistic of zero.  Removing the one from the beginning and relabeling the permutation with the numbers $1$ through $n-1$ gives a permutation in $\Av_{n-1}(\Pi)$ that also has an inversion statistic of zero. Thus $\widetilde{I_n^1}(\Pi;q)=\widetilde{I_{n-1}^1}(\Pi;q)=1$.  

Now suppose the strip tableau has a last column of height $k$ for $2 < k < n$.  For all such shapes, the $1$ must be the element at the top of the last column.  The corresponding permutation is then $w_{column}$ 1 $w_{row}$ where $w_{column}=\begin{array}{cccc} c_1&c_2&\cdots&c_{k-1} \end{array}$ is the word of the last column read from bottom to the second from the top position and $w_{row}= \begin{array}{cccc} r_1&r_2&\cdots&r_{n-k} \end{array}$ is the word of the row read from right to left.  Since the permutation avoids $312$, $r_1 > c_{k-1}$ thus $c_{k-1}=2$.  Since 1 creates an inversion with all of the elements in $w_{column}$, the inversion statistic for the permutation is $\inv(w_{column}) + \inv(w_{row}) + (k-1) + (\text{the \# of inversions formed between the elements in $w_{column}$ and $w_{row}$}).$  Now remove the 1 from this permutation and then relabel the elements from $1$ to $n-1$ in the same relative order.  

If $w_{column} = \begin{array}{cccccc} c_1&c_2&\cdots&c_{k-3}&3&2 \end{array}$ then after removing the $1$ and relabeling we again obtain a permutation that avoids $\Pi$ and thus is a reading word for a strip tableau of size $n-1$ with a single row followed by a last column of height $k-1$.  The inversion statistic for this new tableau is then $\inv(w_{column}) + \inv(w_{row}) +(\text{the \# of inversions formed between the elements in $w_{column}$ and $w_{row}$}).$  Thus the change in the inversion statistic is $k-1$.

If $w_{column}=\begin{array}{ccccc} c_1&c_2&\cdots&c_{k-2}&2 \end{array}$ in the original permutation and $w_{row}=\begin{array}{cccc} 3&r_2&\cdots&r_{n-k} \end{array},$ then the new permutation obtained by removing $1$ is $c_1c_2\cdots c_{k-2}$ 2 3 $ r_2\cdots r_{n-k}$ \ignore{$\begin{array}{ccccccccc} c_1&c_2&\cdots&c_{k-2}&2&3&r_2&\cdots&r_{n-k} \end{array}$} which contains a 312 pattern with $\begin{array}{ccc} c_{k-2}&2&3 \end{array}$ (before relabeling).  We create a new permutation of length $n-1$ by swapping the position of $2$ and $3$.  This new permutation is $c_1c_2\cdots c_{k-2}$ 3 2 $r_2r_3\cdots r_{n-k}$ and avoids the given $\Pi$.  This permutation corresponds to a strip tableau of a single row followed by a column of height $k$.  The inversion statistic for this new permutation is then $\inv(w_{column}) + \inv(w_{row}) + (\text{the \# of inversions formed between the elements in $w_{column}$ and $w_{row}$}) + 1.$  Thus the change in the inversion statistic is $k$.

Therefore, considering that $\widetilde{I_{n}^{n}}(\Pi;q)=q^{\binom{n}{2}}$ and $\widetilde{I_{n}^{2}}(\Pi;q)=q$, we have

\begin{equation}
I_n(\Pi;q) = 1 + q + q^{\binom{n}{2}}+\sum_{k=3}^{n-1}{q^{k-1}\widetilde{I_{n-1}^{k-1}}(\Pi;q) + q^k \widetilde{I_{n-1}^k}(\Pi;q)}.
\label{class9}
\end{equation}
with initial conditions $\widetilde{I_{1}^{1}}(\Pi;q)=\widetilde{I_{2}^{1}}(\Pi;q)=1$ and $\widetilde{I_{2}^{2}}(\Pi;q)=q.$

\subsection*{Class 10 - $\Pi=\{321,312, 132\}$}\label{Class10}

Since permutations in $\Av_n(\Pi)$ avoid the $132$ pattern, they correspond to Fibonacci tableaux for which the topmost entries in any two adjacent columns of height two increase from left to right.  In addition, if a column of height two is followed by a column of height one, the top entries in these two columns also increase from left to right.

If a tableau with these conditions begins with a column of height one, the number in this column must be $n$.  In the corresponding permutation, $n$ is then the last element and thus removing this $n$ gives a permutation corresponding to a permutation of length $n-1$ that avoids $\Pi$ and has the same inversion statistic.  Therefore $I_n^1(\Pi;q) = I_{n-1}(\Pi;q)$.  

If a tableau with these conditions begins with a column of height two, let $k$ be such that the tableau begins with $k$ columns of height two followed by a column of height one.  Because of the condition that the topmost elements of adjacent columns of height two and a column of height one following them must increase left to right, we know that the first $k+1$ columns of the tableau look like
\[
T_{2k+1} = \begin{array}{ccccc}
a_1&a_2&\cdots&a_k& \\
n&{n-1}&\cdots&{n-k+1}&{n-k}
\end{array}\]
where $a_1 < a_2 < \cdots < a_k$.

These first $k+1$ columns are followed by a Fibonacci tableau of size $n-2k-1$, call it $T_{n-2k-1}$, that corresponds to a permutation avoiding $\Pi$ (i.e. $T_{n-2k-1}$ also has the condition about the topmost element in adjacent columns of height two).  The word of the original tableau is then $w_c(T_{n-2k-1}) w_c(T_{2k+1})$.  The number of inversions in $w_c(T_{2k+1})$ is $q^{2(1+2+\cdots+k)} = q^{k(k+1)}$.  The number of inversions that the elements in $T_{2k+1}$ make with $T_{n-2k-1}$ is 

$(n-2k-1) - (a_1-1) + (n-2k-1)- (a_2 - 2) + \cdots + (n-2k-1) - (a_k - k) = k(n-2k-1) - (a_1 + a_2 +\cdots + a_k) + (1+2+\cdots + k)= k(n-2k-1) - (a_1 + a_2 + \cdots + a_k) + \frac{k(k+1)}{2}.$
  Thus 

\begin{align*}
I_n^{2}(\Pi;q) &= \sum_{k=1}^{\lfloor{\frac{n}{2}}\rfloor} \sum_{0 < a_1 < \cdots < a_k < n} q^{k(k+1)}q^{k(n-2k-1) - (a_1 + a_2 + \cdots + a_k) + \frac{k(k+1)}{2}} I_{n-2k-1}(\Pi;q)\\
&= \sum_{k=1}^{\lfloor{\frac{n}{2}}\rfloor} \sum_{0 < a_1 < \cdots < a_k < n} q^{k(n-2k-1) + \frac{3k(k+1)}{2} - (a_1 + a_2 + \cdots + a_k)} I_{n-2k-1}(\Pi;q).
\end{align*}  

and 
\begin{equation}
I_{n}(\Pi;q)=I_{n-1}(\Pi;q)+\sum_{k=1}^{\lfloor{\frac{n}{2}}\rfloor} \sum_{0 < a_1 < \cdots < a_k < n} q^{k(n-2k-1) + \frac{3k(k+1)}{2} - (a_1 + a_2 + \cdots + a_k)} I_{n-2k-1}(\Pi;q)
\label{class10}
\end{equation}

\subsection*{Class 11 - $\Pi=\{321,312,123\}$}\label{Class11}

Since permutations in $\Av_n(\Pi)$ must avoid $123$ they correspond to standard Fibonacci tableaux that have either no columns of height one, a single column of height one that occurs in the rightmost or leftmost position, or two columns of height one that occur in the rightmost and leftmost position.  Thus if $n$ is even, we will only obtain standard Fibonacci tableaux with all columns of height two or with two columns of height one in the rightmost and the leftmost column positions and if $n$ is odd we will obtain standard Fibonacci tableaux with a single column of height one in either the rightmost position or the leftmost position.

To determine $I_n^1(\Pi;q)$, suppose we have a tableau with the above restrictions that begins with a column of height one.  This column must contain $n$, thus in the reading word for this tableau $n$ is in the last position.  Then $n$ does not contribute to the inversion statistic so removing this $n$ does not change the inversion statistic.  For $n \geq  3$, the corresponding tableau then begins with a column of height two and satisfies the restrictions given by $\Pi$.  Thus $I_n^1(\Pi;q) = I_{n-1}^{2}(\Pi;q)$, $n \geq 3$.  

Suppose we have a tableau with the above restrictions that begins with a column of height two.  Then the first column looks like $\begin{array}{c} k\\n \end{array}$ so the corresponding word ends with $n$ $k$.  Removing this $n$ and $k$ and relabeling with the numbers $1$ through $n-2$ gives a permutation that corresponds to a Fibonacci tableau with the above restrictions that again begins with a column of height two (for $n \geq 4$) and the inversion statistic has changed by $k$.  Since we have $n-1$ choices for $k$, we obtain \begin{equation} I_n^{2}(\Pi;q) = (q + q^2 + \cdots + q^{n-1}) I_{n-2}^{2}(\Pi;q) \end{equation} for $n \geq 4$ with $I_{2}^{2}(\Pi;q)=q$ and $I_{3}^{2}(\Pi;q)=q+q^2$. 

Hence, 
\begin{equation}
I_n(\Pi;q) = I_{n-1}^{2}(\Pi;q) + (q + q^2 + \cdots + q^{n-1}) I_{n-2}^{2}(\Pi;q).
\label{class11}
\end{equation}

\subsection*{Class 12 - $\Pi=\{321,213, 132\}$}\label{Class12}

Determining the inversion polynomial for the permutations that avoid this subset of consecutive patterns remains an open question.  Interestingly, this class was the last subset of three consecutive patterns to be enumerated and was done by Kitaev and Mansour \cite{KiM} in 2005, several years after the original enumerative work done by Kitaev \cite{Kit} in 2003.  The authors are currently working on determing this last remaining inversion polynomial.

\subsection*{Class 13 - $\Pi = \{321,312\}$}\label{Class13}

As  mentioned previously, the set of permutations that avoid this $\Pi$ is in one to one correspondence with the set of Fibonacci tableaux and the inversion polynomial for Fibonacci tableaux was given by Cameron and Killpatrick \cite{CaK} in 2013.  In addition, Cameron and Killpatrick proved that the inversion polynomial on this set is symmetric and conjectured to be log concave.

\section{Further Research}

There are many interesting directions for further research in this area.  The authors are working to determine the inversion polynomials for the remaining subset $\Pi$ of size three and the remaining subsets $\Pi$ of size two.  Because these sets have fewer restrictions on the strip tableaux, the recursions are more difficult to determine.  In addition, the authors are looking at subsets $\Pi$ of consecutive patterns of length 4 as well as other subsets of generalized patterns that are not all consecutive patterns.  Since the inversion polynomial on the permutations that avoid the subset in Class 13 is symmetric and conjectured to be log concave, one might wonder which of the other inversion polynomials are symmetric, unimodal or log concave.

\newpage

\begin{longtable}{|l|l|C{4in}|}
\caption[Inversion Polynomials for Permutations Avoiding Consecutive Patterns]{Inversion Polynomials for Permutations Avoiding Consecutive Patterns} \label{bigtable} \\
\hline
\hline
\multicolumn{1}{|c|}{\bf Class} & \multicolumn{1}{|c|}{$\mathbf \Pi$} & \multicolumn{1}{|c|}{$\mathbf I_n(\Pi;q)$} \ignore{& Equation}\\\hline
\endfirsthead

\multicolumn{3}{c}%
{ \tablename\ \thetable{} -- continued from previous page} \\
\hline \multicolumn{1}{|c|}{\bf Class} &
\multicolumn{1}{c|}{$\mathbf \Pi$} &
\multicolumn{1}{c|}{$\mathbf I_n(\Pi;q)$} \\ \hline 
\endhead

\hline \multicolumn{3}{|r|}{{Continued on next page}} \\ \hline
\endfoot

\hline \hline
\endlastfoot

\multirow{2}{*}{1} & $321,312, 132,123$& $
\begin{cases}q^{k(k+1)} + q^{k^2},&\text{if $n=2k+1$}\cr
q^{k(k-1)} + q^{k^2},&\text{if $n=2k$}
\end{cases}
$ \\\cline{2-3}
& $321,231,213,123$& $\begin{cases}q^{\binom{n}{2}-k(k+1)} + q^{\binom{n}{2}-k^2},&\text{if $n=2k+1$}\cr
q^{\binom{n}{2}-k(k-1)} + q^{\binom{n}{2}-k^2},&\text{if $n=2k$}
\end{cases}$  \\\hline\hline

\multirow{4}{*}{2} & $321,312, 231, 132$ & $1+q$ \\\cline{2-2}
& $321, 312, 231, 213$  &     \\\cline{2-3}
&$312, 213,132,123$ &$q^{\binom{n}{2}}+q^{\binom{n}{2}-1}$ \\\cline{2-2}
& $231,213,132,123$ &  \\\hline\hline

3&$312, 213, 231, 132$ & $\begin{cases} 1+q^{\binom{n}{2}}, & n=3 \\ 0, & n>3 \end{cases}$  \\\hline\hline

\multirow{2}{*}{4} & $321,312,213,123$ & $q+q^2$  \\\cline{2-3}
& $321,231,132,123$  & $q^{\binom{n}{2}-1}+q^{\binom{n}{2}-2}$   \\\hline\hline

\multirow{4}{*}{5} & $321,312,213,132$ & $1+q^2+\cdots+q^{n-1}$  \\\cline{2-2}
& $321,231,213,132$  &    \\\cline{2-3}
& $312,231,213,123$  & $q^{\binom{n}{2}}+q^{\binom{n}{2}-2}+\cdots+q^{\binom{n}{2}-n+1}$   \\\cline{2-2}
& $312,231,132,123$&    \\\hline\hline

\multirow{2}{*}{6} & $321, 312,231,123$ &   $\begin{cases}2q^kC_k(q) & \text{if $n=2k+1$}\\ q^{k-1}C_{k-1}(q)+q^kC_k(q) & \text{if $n=2k$}   \end{cases}$ \ignore{(\ref{class6odd}), (\ref{class6even})} \\\cline{2-3}
& $321, 213,132,123$ &  $\begin{cases}2q^{\binom{n}{2}-k}C_k(1/q) & \text{if $n=2k+1$}\\ q^{\binom{n}{2}-k+1}C_{k-1}(1/q)+q^{\binom{n}{2}-k}C_k(1/q) & \text{if $n=2k$}   \end{cases}$  \\\hline\hline

\multirow{2}{*}{7} & $321,312,231$ & $\displaystyle I_{n-1}(\Pi;q) + \sum_{k=1}^{\lfloor n/2 \rfloor}{q^{k-1} C_{k-1}(q)\cdot I_{n-2k+1}(\Pi;q)}$   \\\cline{2-3}
& $213,132,123$ &  $\displaystyle q^{\binom{n}{2}}I_{n-1}(\Pi;1/q) + \sum_{k=1}^{\lfloor n/2 \rfloor}{q^{\binom{n}{2}-k+1} C_{k-1}(1/q)\cdot I_{n-2k+1}(\Pi;1/q)}$  \\\hline\hline

\multirow{4}{*}{8} & $321,312,213$ & $1+q+\cdots+q^{n-1}$  \\\cline{2-2}
& $321, 231,132$  &    \\\cline{2-3}
&  $312, 213, 123$ & $q^{\binom{n}{2}}+q^{\binom{n}{2}-1}+\cdots+q^{\binom{n}{2}-n+1}$    \\\cline{2-2}
& $231, 132, 123$ &   \\\hline\hline

\multirow{4}{*}{9} & $312,231,132$ &  $\displaystyle 1 + \sum_{k=2}^n{q^{k-1}\widetilde{I_{n-1}^{k-1}}(\Pi;q) + q^k \widetilde{I_{n-1}^k}(\Pi;q)}$ \\\cline{2-2}
&  $312, 231, 213$ &    \\\cline{2-3}
& $231,213,132$  &  $\displaystyle q^{\binom{n}{2}} + \sum_{k=2}^n{q^{\binom{n}{2}-k+1}\widetilde{I_{n-1}^{k-1}}(\Pi;1/q) + q^{\binom{n}{2}-k} \widetilde{I_{n-1}^k}(\Pi;1/q)}$  \\\cline{2-2}
& $312, 213, 132$ &   \\\hline\hline

\multirow{4}{*}{10} & $321,312,132$ & $\displaystyle I_{n-1}(\Pi;q)+\sum_{k=1}^{\lfloor{\frac{n}{2}}\rfloor} \sum_{0 < a_1 < \cdots < a_k < n} q^{k(n-2k-1) + \frac{3k(k+1)}{2} - (a_1 + a_2 + \cdots + a_k)} I_{n-2k-1}(\Pi;q)$\\\cline{2-2}
&  $321,231,213$ &    \\\cline{2-3}
 & $231, 213, 123$  &  $\displaystyle q^{\binom{n}{2}}I_{n-1}(\Pi;q^{-1})+\sum_{k=1}^{\lfloor{\frac{n}{2}}\rfloor} \sum_{0 < a_1 < \cdots < a_k < n} q^{ \binom{n}{2}-   k(n-2k-1) - \frac{3k(k+1)}{2} + (a_1 + a_2 + \cdots + a_k)} I_{n-2k-1}(\Pi;q^{-1})$   \\\cline{2-2}
 & $312, 132, 123$ &   \\\hline\hline

\multirow{4}{*}{11} & $321,312,123$ & $I_{n-1}^{2}(\Pi;q) + (q + q^2 + \cdots + q^{n-1}) I_{n-2}^{2}(\Pi;q)$ \\\cline{2-2}
& $321,231,123$  &   \\\cline{2-3}
&  $321,213, 123$ &  $q^{\binom{n}{2}} I_{n-1}^{2}(\Pi;1/q) + \left(q^{\binom{n}{2}-1} + q^{\binom{n}{2}-2} + \cdots + q^{\binom{n}{2}-n+1}\right) I_{n-2}^{2}(\Pi;1/q)$  \\\cline{2-2}
& $321, 132, 123$ &   \\\hline\hline

\multirow{2}{*}{12} & $321,213,132$ &   (\ref{class12}) \\\cline{2-2}
& $312,231,123$ &    \\\hline\hline

\ignore{BEGINIGNORE
\multirow{4}{*}{13} & $321,213$ &   (\ref{class13})\\\cline{2-3}
&   &    \\\cline{2-3}
&   &    \\\cline{2-3}
& &   \\\hline\hline

\multirow{2}{*}{14} & $312,132$ &   (\ref{class14}) \\\cline{2-3}
& $231,213$ &    \\\hline\hline

\multirow{2}{*}{15} & $312,231$ &   (\ref{class15}) \\\cline{2-3}
& $213,132$ &    \\\hline\hline

16 & $321,123$ &   (\ref{class16}) \\\hline\hline
ENDIGNORE}

\multirow{4}{*}{13} & $321,312$ & $I_{n-1}(\Pi;q) + (q + q^2 + \cdots + q^{n-1})I_{n-2}(\Pi;q)$   \\\cline{2-2}
&  $321,231$ &  \\\cline{2-3}
& $213,123$  &  $q^{\binom{n}{2}}I_{n-1}(\Pi;1/q) + \left(q^{\binom{n}{2}-1} + q^{\binom{n}{2}-2} + \cdots + q^{\binom{n}{2}-{n-1}}\right)I_{n-2}(\Pi;1/q)$  \\\cline{2-2}
& $132,123$ &   \\\hline\hline

\ignore{BEGINIGNORE
\multirow{2}{*}{18} & $312,213$ &   (\ref{class18}) \\\cline{2-3}
& $231,132$ &   \\\hline
ENDIGNORE}

\end{longtable}

\end{document}